\documentclass[a4paper,12pt]{article} 
\usepackage[utf8]{inputenc}
\usepackage[T1]{fontenc}
\usepackage{amsmath,amssymb,amsthm}
\usepackage[]{authblk}
\usepackage{bbm}
\usepackage{color}
\usepackage{bbm}
\usepackage{verbatim}
\usepackage[colorlinks]{hyperref}
\usepackage[capitalise]{cleveref}

\usepackage{pgf,tikz}
\usetikzlibrary{arrows}

\newtheorem{thm}{Theorem}
\crefname{thm}{Theorem}{Theorems}

\crefname{cor}{Corollary}{Corollaries}

\crefname{lem}{Lemma}{Lemmas}
\newtheorem{prop}[thm]{Proposition}
\crefname{prop}{Proposition}{Propositions}

\crefname{conj}{Conjecture}{Conjectures}

\crefname{ques}{Question}{Questions}
\theoremstyle{definition}

\crefname{defn}{Definition}{Definitions}

\crefname{rem}{Remark}{Remarks}

\crefname{ex}{Example}{Examples}

\crefname{obs}{Observation}{Observations}

\crefname{claim}{Claim}{Claims}

\crefname{ass}{Assumption}{Assumptions}

\numberwithin{thm}{section}

\newcommand{\cE}{\ensuremath{\mathcal E}}

\newcommand{\cheeger}{\mathrm{Ch}}
\newcommand{\smooth}[1]{\mathrm{C}^{#1}}

\newcommand{\diri}{\mathcal{E}}

\renewcommand{\a}{{\ensuremath{\alpha}}}
\let\oldd\d
\renewcommand{\d}{{\ensuremath{\delta}}}

\newcommand{\E}{{\mathcal{E}}}

\let\oldk\k
\renewcommand{\k}{{\ensuremath{\kappa}}}

\let\oldl\l
\renewcommand{\l}{{\ensuremath{\lambda}}}
\let\oldL\L
\renewcommand{\L}{{\ensuremath{\Lambda}}}

\let\oldo\o
\renewcommand{\o}{{\ensuremath{\omega}}}
\let\oldO\O
\renewcommand{\O}{{\ensuremath{\Omega}}}

\let\oldr\r
\renewcommand{\r}{{\ensuremath{\rho}}}

\let\oldt\t
\renewcommand{\t}{{\ensuremath{\tau}}}
\let\oldu\u
\renewcommand{\u}{{\ensuremath{\upsilon}}}

\renewcommand{\leq}{\leqslant}
\renewcommand{\geq}{\geqslant}
\renewcommand{\le}{\leqslant}

\renewcommand{\to}{\rightarrow}

\begin{document}
\title{Nonlinear Dirichlet forms, energy spaces, and calculus rules}
\author{Giovanni Brigati\thanks{\textsf{brigati@ceremade.dauphine.fr}} }
\affil{CEREMADE, CNRS, Universit\'e Paris-Dauphine, PSL University\protect\\75016 Paris, France}
\date{\vspace{-0.25cm}\today}
\maketitle
\vspace{-0.75cm}

\begin{center}
    \textit{To my late grandfather, Roberto, with love and gratitude.}
\end{center}

\begin{abstract}
We review recent contributions on nonlinear Dirichlet forms. Then,
we specialise to the case of $2-$homogeneous and local forms. Inspired by the theory of Finsler manifolds and metric measure spaces, we establish new properties of such nonlinear Dirichlet forms, which are reminiscent of differential calculus formulae.
\end{abstract}


\section{Introduction}\label{sec:intro}

The theory of (bilinear) Dirichlet forms \cite{bouleau1991dirichlet,fukushima2011dirichlet,ma2012introduction} is a rich topic at the interface between analysis and probability, in connection with Markov semigroups (see Theorem \ref{thm:thm1} below). Dirichlet forms were introduced in \cite{beurling1959dirichlet} as a class of bilinear/quadratic forms generalising the standard Dirichlet energy 
\begin{equation}
    \label{eq:eq1}
    \E(u) = \begin{cases}
         \int_{\mathbb R^d} |Du(x)|^2 \, dx, \qquad \text{if} \, u \in \mathrm{W}^{1,2}(\mathbb R^d); \\
        +\infty, \qquad \text{otherwise}.
    \end{cases}
\end{equation}
An abstract calculus (called $\Gamma$-calculus) has been developed in \cite{bakry1985diffusions} to capture the hypercontractivity and decay properties of linear Markov semigroups induced by quadratic Dirichlet forms, establishing connections between linear diffusion processes, Riemannian geometry, and functional inequalities \cite{bakry2014analysis,dolbeault2008bakry,arnold2001convex}.
In particular, on Riemannian manifolds, the $\Gamma$-calculus applied to the $\mathrm{W}^{1,2}$-seminorm (which is a quadratic Dirichlet form) links the long-time behaviour of the heat flow (the associated Markov semigroup) with the Bochner identity \cite{DOCARMORIEMGEOM}, via an inequality called \emph{Bakry--Emery condition}. In case of a lower bound on the \emph{Ricci curvature}, the Bakry--Emery condition is satisfied. This turns out to be an equivalence, yielding a definition of Ricci lower bounds in pure terms of Dirichlet forms \cite{wang2011equivalent}.
A third equivalent notion of Ricci lower bounds was given via optimal transport by Lott, Villani, and Sturm \cite{lott2009ricci,sturm2006geometryI,sturm2006geometryII}.

This last definition makes sense even for metric measure spaces. However, unlike Riemannian manifolds, the analogous of \eqref{eq:eq1} on metric measure spaces, called \emph{Cheeger's energy} \cite{ambrosio2013density,ambrosio2014calculus,ambrosio2014metric}, is not a quadratic Dirichlet form in general, and the metric heat flow is nonlinear. So, $\Gamma$-calculus and a Bakry--Emery condition could be established only in metric measure spaces whose Cheeger's energy is quadratic (\emph{infinitesimally Hilbertian spaces})  \cite{ambrosio2015bakry}. In this class, Ambrosio, Gigli, and Savar\'e could recover the equivalence between the Bakry--Emery condition and the Ricci lower bounds of \cite{lott2009ricci,sturm2006geometryI}. 
Then, the study of geometry and functional inequalities on RCD spaces (i.e. infinitesimally Hilbertian spaces satisfying a Ricci lower bound) flourished, becoming a main subject in the last ten years \cite{gigli2015differential,gigli2018nonsmooth}. 

Much less is known in case the Cheeger's energy is nonquadratic. At the moment, Ricci lower bounds in this case are given only via the Lott--Villani--Sturm approach. The Cheeger's energy being nonquadratic is not exceptional, as it is the case in all (non-Riemman) Finsler geometries. Finsler structures will be a source of inspiration for this paper, as Riemmanian manifolds are model examples for the RCD case. 
We detail hereby the state of the art. 

\begin{itemize}
    \item Ricci lower bounds on Finsler manifolds appear in \cite{ohta2009heat,ohta2014bochner,ohta2021comparison}. The equivalence between an \emph{intrinsic definition}, and the Lott--Villani--Sturm approach of \cite{lott2009ricci,sturm2006geometryI} holds true. Moreover, a suitable equivalent Bakry--Emery condition has been found. This last condition looks very different from the standard one, as it is a comparison estimate between the nonlinear Finsler Laplacian and a linearisation of it. 
    \item A definition of \emph{nonlinear Dirichlet form} has been given in \cite{cipriani2003nonlinear}. Properties of nonlinear Dirichlet forms have been studied in \cite{claus2021energy,claus2021nonlinear} and \cite{brigati2022normal,brigati2022kinetic}. No analogous of the Bakry--Emery $\Gamma$-calculus is available in this context up to the best of our knowledge.
    \item Nonquadratic Cheeger's energies belong to the class of nonlinear Dirichlet forms of \cite{cipriani2003nonlinear}. The converse, i.e.~a representation theorem of abstract nonlinear Dirichlet forms as metric Cheeger's energies, is missing in general, but available in the quadratic case \cite{ambrosio2015bakry}.
\end{itemize}
Motivated by the last point, our long-term goal is to strengthen the link between nonlinear Dirichlet forms and nonquadratic Cheeger's energies. In this note, we study $2$-homogeneous, and local nonlinear Dirichlet forms. We analyse the associated non-Hilbertian energy space, and we recover some abstract calculus rules, which are reminiscent of \emph{concrete} calculus in metric measure spaces \cite{gigli2015differential,ambrosio2014calculus}, capturing the analogies as much as we can.

The work is organised as follows. Section \ref{sec:sec2} contains a presentation of the main notion involved in the paper. Section \ref{sec:sec3} specialises to $2$-homogeneous and local nonlinear Dirichlet forms. Property of the associated energy spaces are collected in Section \ref{sec:sec4}. Sections \ref{sec:sec5}-\ref{sec:sec6} contain first-order, and second-order calculus rules, respectively. Finally, in Section \ref{sec:sec7} we list some desirable results which are still missing in the theory. 
 
\section{Definitions and tools}\label{sec:sec2}

\subsection{Quadratic Dirichlet forms}

Let $(X,\mathcal{F},m)$ be a $\sigma$-finite measure space, such that $\mathcal{F} $ is in a bi-measurable correspondence with the Borel class on $\mathbb R$.

A (quadratic) Dirichlet form is a quadratic, lower-semicontinuous (l.s.c.) functional $$\E: \mathrm{L}^2(X,m) \to [0,\infty],$$
whose domain
$$D(\E) := \left\{u \in \mathrm{L}^2(X,m) \, : \, \E(u) \neq \infty\right\}$$
is a dense subspace of $\mathrm{L}^2(X,m)$, and such that
\begin{equation}\label{eq:eq5}
\forall u \in D(\E), \quad \E(0 \vee u \wedge 1) \leq \E(u).
\end{equation}
The symbols $\vee, \, \wedge$ stand for the maximum and minimum operations, respectively.
The quadratic form $\E$ induces an unbounded, positive semi-definite, self-adjoint linear operator $A: D(A) \to \mathrm{L}^2(X,m),$ such that 
$$\forall u \in D(\E), \qquad \E(u) = \int_X \sqrt{A} u \, \sqrt{A}u \, dm,$$
and 
$$\forall u \in D(A), \qquad 2\E(u) = \int_X -u \, Au \, dm.$$
By construction, one has that $t \mapsto T_t:= \mathrm{e}^{At}$ is a linear and continuous semigroup of contractions on $\mathrm{L}^2(X,m)$ \cite{ccinlar2011probability} such that 
$$\partial_t T_t = A \, T_t.$$ 
Finally, one could introduce the bilinear Dirichlet form associated with $\E$ as follows 
\begin{equation}\label{eq:eq6}
\Lambda(u,v) := \int_X \sqrt{A}u \, \sqrt{A}v, \qquad \forall u,v \in D(\E),
\end{equation}
which is called simply \emph{Dirichlet form} in the literature \cite{bouleau1991dirichlet,fukushima2011dirichlet,ma2012introduction}. 

The interest of Dirichlet forms, especially in connection with linear semigroups \cite{fukushima2011dirichlet}, lies in the following.

\begin{thm}\label{thm:thm1}
    Let $\E$ be a quadratic, l.s.c., densely defined, positive semi-definite quadratic form over $\mathrm{L}^2(X,m)$.
    Then, the following are equivalent.
    \begin{enumerate}
        \item $\E$ is a Dirichlet form.
        \item The linear semigroup $(T_t)_t = \mathrm{e}^{At}$ is a continuous and self-adjoint semigroup of contractions in $\mathrm{L}^p(X,m)$ for all $p \in [1,\infty]$:
        \begin{equation}\label{eq:eq3}
\forall t \geq 0, \quad \forall p \in [1,\infty], \qquad \|T_t\|_{\mathrm{L}^p(\mathbb R^d) \to \mathrm{L}^p (\mathbb R^d)} \leq 1.
\end{equation}
        Moreover, $(T_t)_t$ is positivity-preserving: 
        $$\forall u \in \mathrm{L}^2(X,m), \, t \geq 0, \qquad u \geq 0 \, \implies \, T_t u \geq 0.$$
        \item The form $\E$ verifies the normal contraction property \begin{equation}\label{eq:eq4}
    \forall u \in \mathrm{W}^{1,2}(\mathbb R^d), \quad \forall \phi \in \Phi, \qquad \E(\phi(u)) \leq \E(u), 
\end{equation}
where 
$$\Phi := \left\{ \phi : \mathbb R \to \mathbb R \, : \, \phi(0)=0, \, \phi  \text{ is }  \text{$1$-Lipschitz} \right\}.$$.
    \end{enumerate}
\end{thm}

The implication $(2 \implies 1)$ should be understood in the following sense. Given a linear, self-adjoint semigroup $(T_t)_t$, satisfying the hypotheses of condition $2.,$ one could always (uniquely) define an unbounded, linear, positive semi-definite, and self-adjoint operator 
$$A := \lim_{t \to 0} \, \frac{T_t - \mathrm{Id}}{t},$$
then, $\sqrt{A}$ by functional calculus. 
This way, the functional 
$$\E(u) := \begin{cases}
    \int_X \left| \sqrt{A} u \right|^2 \, dm, \qquad \text{if } u \in D(\sqrt{A}), \\
    +\infty, \qquad \text{otherwise}, 
\end{cases}$$
is a Dirichlet form.
Semigroups satisfying condition $2.$ of Theorem \ref{thm:thm1} are usually called linear \emph{Markov semigroups}.

\subsection{Nonlinear Dirichlet forms}

In \cite{cipriani2003nonlinear}, a possible extension of Dirichlet forms to the nonlinear setting has been established.
We adopt this approach, in view of a nonlinear extension of Theorem \ref{thm:thm1}, namely Theorems \ref{thm:thm2} and \ref{th:contraction}.

Let $\E: \mathrm{L}^2(X,m) \to [0,\infty]$ be a convex, l.s.c.~functional with dense domain. Indicate with $\partial \E:  D(\partial \E) \subset \mathrm{L}^2(X,m) \to 2^{\mathrm{L}^2(X,m)}$ its subdifferential operator:
$$\partial \E(u) = \left\{ \xi \in \mathrm{L}^2(X,m) \, : \, \forall z \in \mathrm{L}^2(X,m), \quad \E(z) - \E(u) \geq \int_X \xi \, u \, dm \right\},$$
see \cite{brezis1973operateurs,brezis2011functional}.
Let $(T_t)_{t\geq0}$ be the semigroup of nonlinear operators generated by $-\partial \cE$ via the differential equation
\begin{equation}
\label{gradientflow}
\begin{cases}
    \partial_t T_t \,u  \in - \partial \E( T_t \, u),&\forall t \in (0,\infty), \quad \forall u \in \mathrm{L}^2(X,m),\\
    T_0 \,u = u,&\forall u \in \mathrm{L}^2(X,m). 
\end{cases}    
\end{equation}
\Cref{gradientflow} is well-posed for all $u\in\mathrm{L}^2(X,m)$. Its solution is usually called the gradient flow of $\E$ starting at $u$. See \cite{ambrosio2008gradient,brezis1973operateurs}.

The functional $\E$ is a nonlinear Dirichlet form if the associated semigroup $(T_t)_t$ is a contraction in all $\mathrm{L}^p(X,m)$ spaces \eqref{eq:eq3}, and if, in addition, it is order-preserving: 
\begin{equation}
    \label{eq:eq7}
    \forall u,v \in \mathrm{L}^2(X,m), \, t \geq 0, \qquad u \geq v \, \implies \, T_t u \geq T_t v.
\end{equation}

A semigroup of nonlinear maps satisfying \eqref{eq:eq3}, and \eqref{eq:eq7} is called a \emph{nonlinear Markov semigroup}.
Notice that a quadratic Dirichlet form is a special case of a nonlinear Dirichlet form. Moreover, if $\E$ is quadratic, we have that $\partial \E =A$.

Thanks to \cite{benilan1979quelques}, condition \eqref{eq:eq3} should be verified only for $p=\infty$. After the results in  \cite{barthelemy1996invariance,brezis1973operateurs,brigati2022normal}, conditions \eqref{eq:eq3},\eqref{eq:eq7} can be characterised in terms of invariance of convex sets in $\mathrm{L}^2(X,m)$ for the action of the semigroup $(T_t)_t$, then, equivalently rewritten as contraction properties of the functional $\E$ itself. 

\begin{thm}[\cite{brigati2022normal}]\label{thm:thm2}
Let $\E:\mathrm{L}^2(X,m)\to[0,\infty]$ be a 
l.s.c.\ functional. 
Then, $\mathcal E$ is a nonlinear Dirichlet form if and only if, for all $u,v \in \mathrm{L}^2(X,m)$, and $\a\in[0,\infty)$, $\E$ verifies 
\begin{align}
\label{eq:minmax}
\cE(u\vee v)+\cE(u\wedge v)\le{}&\cE(u)+\cE(v),\\
\label{eq:HK}
\cE(H_\a(u,v))+\cE(H_\a(u,v))\le{}&\cE(u)+\cE(v),
\end{align}
with
\begin{equation}
\label{eq:def:H}
H_\a(u,v)(x)=\begin{cases}
v(x)-\a&u(x)-v(x)<-\a,\\
u(x)&u(x)-v(x)\in[-\a,\a],\\
v(x)+\a&u(x)-v(x)>\a.
\end{cases}
\end{equation}
\end{thm}
The \emph{normal contraction property} \eqref{eq:eq4} was recovered also for the nonlinear setting.
\begin{thm}[\cite{brigati2022normal}]\label{th:contraction}
Let $\E$ be a nonlinear Dirichlet form. Then $\E$ has the normal contraction property \eqref{eq:eq4} if and only if
\begin{equation}
    \label{eq:sym}
    \E(-f) \leq \E(f)\quad\forall f \in \mathrm{L}^2(X,m).
\end{equation}
\end{thm}
By homogeneity, condition \eqref{eq:sym} is trivial in the quadratic setting. 

\subsection{Finsler and metric Sobolev spaces}

Let $\mathcal{M}$ be a smooth, closed, boundary-free manifold. Let $m$ be a Borel, $\sigma$-finite measure on $\mathcal{M}$. A Finsler structure on $\mathcal{M}$ is a smooth map $x\mapsto |\cdot|_x$ associating to each point $x \in \mathcal{M}$ a norm on $T_x\mathcal{M}.$ Let $|\cdot|^\star_x$ be the dual norm, and consider the map $F : T\mathcal{M} \to 2^{T^\star \mathcal{M}}$ defined as follows:
$$F(x,\xi) = \left\{ (x,\zeta) \, : \, \zeta \in T^\star_x\mathcal{M}, \,   |\zeta|^\star_x = |\xi|_x, \quad \langle \zeta ,\xi \rangle_x = |\xi|_x^2  
\right\},$$
where $\langle \cdot, \cdot \rangle_x$ is the duality product between $T^\star_x\mathcal{M}$ and $T_x\mathcal{M}.$
If $u : \mathcal{M} \to \mathbb R$ is a regular function, then $Du$ is a section of $T^\star\mathcal{M}.$

We have that the Sobolev seminorm 
\begin{equation}
    \label{eq:eq8}
    \E(u) =\begin{cases} \int_\mathcal{M} \left( |Du(x)|_x^\star\right)^2 \, dm = \int_{\mathcal{M}} \left| F^{-1}(Du)\right|^2_x \, dm, \quad \text{if } u \in \mathrm{W}^{1,2}(\mathcal{M},dm),  \\
    \infty, \qquad \text{otherwise},
    \end{cases}
\end{equation}
is a nonlinear Dirichlet form. 
Moreover, $\E$ is a quadratic form if and only if $F$ is a linear map, if and only if $|\cdot|_x$ is induced by a scalar product for all $x \in \mathcal{M}$, so that $\mathcal{M}$ is a Riemmanian manifold.
One can compute 
$$\partial \E(u) = - 2\nabla \cdot F^{-1} (Du),$$
which plays the role of the Laplacian on $(\mathcal{M},|\cdot|,m)$, but it is a nonlinear (possibly multivoque) operator.
If $u \in D(\partial\E)$, and $v \in \mathrm{L}^2(\mathcal{M},m)$, then the scalar product $$\int_{\mathcal{M}} - \nabla \cdot F^{-1} (Du) \, v \, dm$$
makes sense. As Gigli observes in \cite{gigli2015differential}, one cannot hope to find an adjoint operator $S$ of $- \nabla \cdot F^{-1} (Du)$ such that 
$$\int_{\mathcal{M}} - \nabla \cdot F^{-1} (Du) \, v \, dm = \int_{\mathcal{M}} S(v) \, u \, dm,$$
as the right-hand-side is linear in $u$, while the left-hand-side is not. The best one could do is moving only one derivative on $v$ and finding 
$$\int_{\mathcal{M}} - \nabla \cdot F^{-1} (Du) \, v \,  dm = \int_{\mathcal{M} } \langle F^{-1} (Du) , Dv \rangle_x \, dm .$$ 
The form 
\begin{equation}
    \label{eq:eq9}
    \Lambda(u,v) := \int_{\mathcal{M} } \langle F^{-1} (Du) , Dv \rangle_x \, dm
\end{equation}
is multivoque, defined for $u,v \in D(\E)$. The two arguments play different roles in $\Lambda,$ but there is no hierarchy in their regularity. 
We define also the maximal and minimal sections of $\Lambda,$ given by 
\begin{eqnarray}
    \label{eq:eq10}
    &\Lambda^{-} (u,v) := \int_{\mathcal{M}} \inf_{g \in F^{-1}(Du)} \langle g, Dv \rangle _x \, dm, \qquad \forall u,v \in D(\E), \\
    \label{eq:eq11}
    &\Lambda^{+} (u,v) := \int_{\mathcal{M}} \sup_{g \in F^{-1}(Du)} \langle g, Dv \rangle _x \, dm, \qquad \forall u,v \in D(\E).
\end{eqnarray}
Finally, we have \cite{gigli2015differential}
\begin{equation}\label{eq:eq12}
\forall u,v \in D(\E), \qquad \Lambda^{\pm}(u,v) = \lim_{\sigma \to 0^{\pm}} \frac{\E(u+\sigma v) - \E(u)}{\sigma}.
\end{equation}
We remark that the last formula holds even at a pointwise level in the context of Finsler manifolds.

We conclude the section by recalling the definition of the Cheeger's energy, which extends \eqref{eq:eq1} and \eqref{eq:eq8} to metric measure spaces \cite{ambrosio2014calculus}.

Let $(X,\tau)$ be a topological space such that it is homeomorphic to a complete and separable metric space. Then $(X,\tau)$ is called a Polish space. Let $(X,\tau)$ be a Polish space equipped with a function $d: X \times X \to [0,+\infty]$ such that
\begin{itemize}
\item
$d$ is an extended distance on $X$;
\item
$d$ is $\tau-$l.s.c.;
\item
for all sequences $(x_n)_n \subset X$ such that $d(x_n,x) \to 0$, for an element $x \in X$, we have $x_n \to  x$ in $\tau$.
\item
the extended metric space $(X,d)$ is complete. 
\end{itemize}
Let $m$ be a Borel, $\sigma-$finite measure on $(X,d,\tau)$ such that 
\begin{equation}\label{amb}
m(B(x,r)) \leq \exp(Cr^2),
\end{equation}
for a uniform constant $C$.
Let $f \in \mathrm{Lip}_b(X)$ and let $x \in X.$
The maximum local slope of $f$ at $x$ is given by
$$|Df|(x) := \limsup_{y \to x} \dfrac{|f(y)-f(x)|}{d(x,y)}.$$
Let $\bar{\mathrm{Ch}}: \mathrm{L}^2(X,m) \to [0,+\infty]$ be defined via 
$$\bar{\mathrm{Ch}}(u) = \dfrac{1}{2} \int_X |Du|^2 dm \qquad \text{ if $u \in \mathrm{Lip}_b(X)$},$$
and $+\infty$ otherwise. 

The Cheeger's energy of $(X,d,\tau,m)$ is the functional $\mathrm{Ch}$ defined via
$$\mathrm{Ch} = sc^- \bar{\mathrm{Ch}},$$
where $sc^-$ is the l.s.c.~envelope in the $\mathrm{L}^2$-topology.

We have that $D(\mathrm{Ch})$ is a vector space, known as the metric Sobolev space and indicated in the literature with the symbol $\mathrm{W}^{1,2}(X,d,m)$. In general, the Cheeger's energy $\mathrm{Ch}$ is not quadratic and $\mathrm{W}^{1,2}(X,d,m)$ is not a Hilbert space. Let $u \in D(\partial \mathrm{Ch})$. Then,
$$-\Delta_{d,m}(u) := \partial^0 \cheeger(u),$$
is called the metric Laplacian of $u$, where $\partial^0$ indicates the element of minimal norm in the subdifferential. Our analysis of nonlinear Dirichlet forms is motivated by the following.

\begin{thm}[\cite{ambrosio2014calculus}]
    The Cheeger's energy $\mathrm{Ch}$ is a nonlinear Dirichlet form in the sense of Cipriani and Grillo \cite{cipriani2003nonlinear}.
\end{thm}

For calculus rules in metric spaces we generally refer to \cite{ambrosio2014calculus,gigli2015differential}. 
Precise results will be cited wherever they are needed.

\section{\texorpdfstring{$2$}{}-homogeneous and local functionals}\label{sec:sec3}

Even when the Cheeger's energy is nonquadratic, it is a $2$-homogeneous functional satisfying some locality properties \cite{ambrosio2014calculus}. 
Then, we reduce to the class of nonlinear Dirichlet forms which are $2$-homogeneous and local (in a sense defined below). 

Let $\diri: \mathrm{L}^2(X,m) \to [0,\infty]$ be a functional. Then, $\diri$ is $2-$homogeneous if 
$$\forall \nu \in \mathbb R, \, \forall u \in \mathrm{L}^2(X,m), \qquad \diri(\nu u) = \nu^2 \diri(u).$$

We say that $\E$ is local if, for all $u,v \in D(\E)$ such that $u$ is constant on the support of $v$, we have
$$\E(u+v) = \E(u) + \E(v).$$ 

Finally, we introduce the symbol $(\cdot,\cdot)$ as a short-hand notation for the $\mathrm{L}^2(X,m)-$scalar product.

\begin{thm}
    \label{thm:thm3}
Let $\E$ be a $2$-homogeneous, convex, l.s.c. functional. Then, the following hold true.

\begin{enumerate}
\item $\E(0) = 0.$
    \item The subdifferential $\partial \diri$ is $1-$homogeneous. Moreover the set $D(\partial\diri)$ is invariant under scalar multiplication. Conversely, if $\diri$ is a nonnegative, convex, l.s.c. functional, such that its subdifferential $\partial \diri$ is $1-$homogeneous and $\diri(0)=0$, then $\diri$ is $2-$homogeneous.
    \item For all $u \in D(A),  \,\forall \xi \in \partial\diri(u)$ we have 
\begin{equation}\label{cdc2}
\int_X u \xi \, dm = 2 \diri(u). 
\end{equation}
\end{enumerate}

\end{thm}

Formula \eqref{cdc2} is reminiscent of the integration by parts 
$$\int_{\mathbb R^d} |Du|^2 \, dx = - \int_{\mathbb R^d} u \, \Delta u \, dx.$$ 
In Finsler geometry, the same holds: 
$$ \E(u) = \int_{\mathcal{M}} \langle F^{-1}(Du), Du \rangle_x \, dm = - \int_{\mathcal{M}} \nabla \cdot F^{-1}(Du) \, u \, dm, \qquad \forall u \in D(\partial \E).$$

\begin{proof}[Proof of Theorem \ref{thm:thm3}]
The first property is trivial:
$\E(0) = 0^2 \E(0) = 0.$

For the second property, suppose that $u \in D(\diri)$ and let $\lambda >0$.
Let $w \in \partial \diri(u)$, so that
$$\diri(v)-\diri(u) \geq (w, v-u),$$
for all $v \in \mathrm{L}^2(X,m)$. 
Then 
$$\diri(\lambda v) -\diri(\lambda u) = \lambda^2(\diri(v)-\diri(u))$$
hence
$$\diri(\lambda v) -\diri(\lambda u)\geq \lambda^2(w,v-u)=(\lambda w, \lambda v - \lambda u).$$
Since $v$ is arbitrary, so it is $\lambda v$.

For the converse, suppose in addition that $\E$ is a $\mathrm{C}^{1,1}$ functional. For a fixed element $u$, it is sufficient to prove that $t \mapsto \diri(tu)$ is $2-$homogeneous. This fact is straightforward, since a real $\mathrm{C}^1$ function which vanishes in $0$, and whose derivative is $1-$homogeneous, is bound to be $2-$homogeneous. In the general case, one can argue by Yosida regularisation \cite{evans1998partial}, provided that the Yosida-regularised of $\E$ is still $2$-homogeneous, as we show below. 

For the third property, we use Yosida regularisation, with the notation of \cite{evans1998partial}. In particular, let $A_\lambda$ be the regularisation \emph{\`a la Yosida} of $\partial \E$, for all $\lambda >0.$ If $\diri$ is $2-$homogenous, then so it is its Yosida-regularised (and vice-versa). Indeed:
\begin{align*}
\diri_\lambda(\mu u) &= \\
                     &= \inf_{v \in \mathrm{L}^2(X,m)} \left\{ \diri(\mu v) + \dfrac{1}{2 \lambda} | \mu u- \mu v|^2\right\} = \\
                     &= \mu^2 \diri_\lambda(u),
\end{align*}
for all $\mu, \lambda > 0.$
As an intermediate result we prove an explicit formula for $\diri_\lambda.$ 
Fix $u\in D(\partial \diri)$.
Let $g: [0, +\infty) \to \mathbb{R}$ be defined via 
$$g(x)=\diri_\lambda(xu)-x^2\diri_\lambda(u).$$
We have that $g$ is $\smooth{1}$ by composition and $g=g'=0$, which reads as:
$$(A_\lambda(xu),u)-2x\diri_\lambda(u))=0,$$
that is \eqref{cdc2} for $\diri_\lambda,$
hence 
$$\diri_\lambda(u) = \dfrac{1}{2} (A_\lambda(u),u).$$
We can prove the direct implication. If one passes in the limit for $\lambda \downarrow 0$, one obtains 
$$\diri(u) = \dfrac{1}{2}(A^0(u),u).$$
Take any other element $\xi \in \partial\diri(u)$. Via $G-$convergence \cite{attouch1984variational}, we find a sequence $(u_n)_n$ s.t. $u_n \to u$ strongly in $\mathrm{L}^2$ and $A_{\frac{1}{n}}(u_n) \to \xi$ strongly in $\mathrm{L}^2$. Without loss of generality choose $(u_n)_n$ such that $|u_n-u| < n^{-2}$. Then
$$\dfrac{1}{2}\left(A_{\frac{1}{n}}(u_n),u_n\right) \to \dfrac{1}{2} (\xi,u).$$
At the same time,  
$$\dfrac{1}{2}\left(A_{\frac{1}{n}}(u_n),u_n\right) = \diri_{\frac{1}{n}}(u_n) \to \diri(u),$$ 
see \cite{evans1998partial,brezis1973operateurs}.

\end{proof}
    
\section{Energy spaces}\label{sec:sec4}

Let $\E$ be a $2$-homogeneous  nonlinear Dirichlet form. 
Then, the space $D(\E)$ is called \emph{Dirichlet space}. Some properties of $D(\E)$ have already been given in \cite{claus2021nonlinear}, but, under our hypotheses, we have a simpler structure with new results. 

\begin{thm}
    Let $\E$ be a $2$-homogeneous nonlinear Dirichlet form. Then, the following hold true.
    \begin{enumerate}
        \item The space $D(\E)$ is a vector space. The functional 
        $$u \mapsto \|u\|^2_{\E} := \|u\|^2_{\mathrm{L}^2(X,m)} + \E(u)$$ is the square of a norm on $D(\E)$. Moreover, $(D(\E),\|u\|_{\E})$ is a Banach space. 
        \item The pair $(D(\E),\|u\|_{\E})$ is a Hilbert space if and only if $\E$ is a quadratic Dirichlet form. 
        \item Lipschitz functions $\phi: \mathbb R \to \mathbb R$ such that $\phi(0)=0$ act on $D(\E)$:
        $$\forall u \in D(\E), \qquad \phi \circ u \in D(\E).$$
        Moreover, $D(\E) \cap \mathrm{L}^\infty(X,m)$ is an algebra. 
        \item $(D(\E),\|u\|_{\E})$ is a dual space \cite{claus2021energy}.
    \end{enumerate}
\end{thm}

\begin{proof}
    Since $D(\E)$ is homogeneous and convex, it is a vector space. The stability under scalar multiplication is encoded in the $2-$homogeneity.
    We shall now prove the triangle inequality for the functional $\|u\|_\diri$, while the homogeneity and the condition $\|u\|_\diri = 0$ iff $u=0$ are clear. For any normed vector space, the triangular inequality is equivalent to the convexity of the closed unit ball. Let $u,v \in D(\diri)$ such that $\|u\|_\diri \vee \|v\|_\diri \leq 1.$
    Then, for any $\lambda \in [0,1]$:
    \begin{align*}
    \|\lambda(u)+(1-\lambda)v\|^2_\diri &\leq \\
    &= |\lambda u +(1-\lambda)v|^2 + \diri(\lambda u+(1-\lambda)v) \leq \\
    &\leq \lambda |u|^2 + (1-\lambda)|v|^2 + \lambda \diri(u) + (1-\lambda)\diri(v) = \\
    &= \lambda(|u|^2+\diri(u))+(1-\lambda)(|v|^2+\diri(v)) \leq \\
    &\leq \lambda + 1-\lambda =1.
    \end{align*}
    Let now $(u_n)_n$ be a Cauchy sequence in $D(\diri)$ with respect to the norm $\|\cdot\|_\diri$. Then $(u_n)_n$ is a Cauchy sequence in $\mathrm{L}^2$. Hence $u_n \rightarrow u$ in $\mathrm{L}^2$, for an element $u$. For $m \in \mathbb{N}$, $\diri(u_n -u_m)$ is bounded in $\mathbb{R}$. Thanks to l.s.c. 
   $$\diri(u - u_m) \leq \liminf_{n \to \infty} \diri(u_n - u_m) \leq C.$$
   Hence, $\diri(u) \leq 2\diri(u-u_m) +\diri(u_m) \leq 4C,$ which implies $u \in D(\diri)$. 
   Finally, 
   $$0 \leq \lim_m \diri(u-u_m) \leq \lim_m \liminf_n \diri(u_n -u_m) \downarrow 0.$$ 
   The first statement is then proved. 
   The second statement follows by definition of $\|u\|_\diri$. 
   For the third statement, if $\phi$ is $1$-Lipschitz, then, Theorem \ref{th:contraction} ensures $u \in D(\E) \, \implies \phi(u) \in D(\E)$. Otherwise, if $\phi$ is $L$-Lipschitz, then $\phi/L$ is $1$-Lipschitz, so 
   $$\E(\phi(u)) = L^2 \, \E(L^{-1}\phi(u)) \leq L^2 \, \E(u).$$
   If $u \in D(\E) \cap \mathrm{L}^\infty$, then $u^2$ is a Lipschitz transformation of $u$. Hence, the function $u^2 \in D(\E).$ The computation $uv = 1/2 \left( (u+v)^2 - u^2 - v^2 \right)$ concludes the proof of the third statement of the theorem. 
\end{proof}

Notice that the third statement of the last theorem replicates a well-known result in the theory of Sobolev spaces. 

\section{First-order calculus}\label{sec:sec5}

Let $\E$ be a $2-$homogeneous nonlinear Dirichlet form. Our goal is to reconstruct an object which plays the role of \eqref{eq:eq10}-\eqref{eq:eq11}, but expressed in pure terms of $\E.$ 
In case $\E$ is quadratic, the natural associated bi-variate object is the bilinear Dirichlet form \eqref{eq:eq6}, and the two limits \eqref{eq:eq10}-\eqref{eq:eq11} coincide. In the general case, a canonical bi-variate object lacks, so we have to take a choice.   

Having \eqref{eq:eq12} in mind, we define 
$$\Lambda^{\pm}(u,v) := \lim_{\sigma \to 0^\pm} \, \frac{\E(u+\sigma v)-\E(u)}{\sigma}, \qquad \forall u,v \in D(\E),$$
as the left and right slopes of $\E$, at $u$, in the direction of $v$. Notice that the definition is well-given, as $\sigma \mapsto \E(u + \sigma v): \mathbb R \to \mathbb R$ is convex. 

\begin{thm}
    \label{thm:thm5} 
    Let $\E$ be a $2$-homogeneous nonlinear Dirichlet form.
    Then, the slopes $\Lambda^\pm$ have the following properties. 
    \begin{enumerate}
        \item For all $u,v \in D(\E)$,  $\Lambda^\pm(u,v)$ are finite, and
        \begin{equation}\label{eq:eq13}
            \left| \Lambda^\pm(u,v) \right| \leq 2 \, \sqrt{\E(u)} \, \sqrt{\E(v)},
        \end{equation}
        which is sharp for $u=v$.
        Moreover,
        $$\Lambda^{-}(u,v) \leq \Lambda^+(u,v), \qquad \forall u,v \in D(\E).$$
        \item For all $u \in D(\E)$, we have $\Lambda^{\pm}(u,u) = 2 \E(u).$
        \item
        For all $u \in D(\E)$, $D^\pm(u,\cdot)$ is positively $1-$homogeneous. Moreover $D^+(u,\cdot)$ is convex, while $D^-(u,\cdot)$ is concave.
        \item
        For all $v \in D(\E)$, we have that $D^\pm(\cdot,v)$ is positively $1-$homogeneous.
        \item
        For all $u,v \in D(\E)$, it holds  $D^+(u,-v)=-D^-(u,v).$ 
        \item If $\E$ is local and $u,v \in D(\E)$ are such that $u$ is constant on $supp(v)$, we have
        $$\Lambda^\pm(u,v)=0.$$
        \end{enumerate}
\end{thm}

\begin{proof}
$1.$\\

    The inequality $\Lambda^- \leq \Lambda^+$ is a consequence of convexity for the map $t \mapsto \E(u+tv).$ Take now any $u,v \in D(\E).$
\begin{align*}
\Lambda^+(u,v) &= \\
		 &=\lim_{h \to 0^+} h^{-1}(\diri(u+hv)-\diri(u)) = \\ &= \lim_h h^{-1}\left(\left(\sqrt{\diri(u+hv)}\right)^2-\left(\sqrt{\diri(u)}\right)^2\right) \leq \\
         &\leq \lim_h h^{-1}\left(\left(\sqrt{\diri(u)}+h\sqrt{\diri(v)}\right)^2-\left(\sqrt{\diri(u)}\right)^2\right) \leq \\
         &\leq 2\sqrt{\diri(u)}\sqrt{\diri(v)} +\lim_h h\diri(v) = 2\sqrt{\diri(u)}\sqrt{\diri(v)}.
\end{align*}
We used the monotonicity of $t \mapsto t^2$ and the fact that $\sqrt{\diri}$ is a seminorm.
With an analogous argument, one can prove that 
$$\Lambda^-(u,v) \geq -2 \sqrt{\diri(u)}\sqrt{\diri(v)}.$$

\noindent $2.$\\

Compute 
$$\Lambda^\pm(u,u) = \lim_{h \to 0^\pm} \frac{\E((1+h)u)-\E(u)}{h} = \lim_{h \to 0^\pm} 2h \, \E(u) + h^2\, \E(u) = 2\, \E(u).$$

\noindent $3.$ and $4.$\\
Let $\lambda >0$. 
\begin{align*}
\Lambda^\pm(u,\lambda v) &=\lim_{h \to 0^\pm}h^{-1}(\diri(u+h\lambda v)-\diri(u)) = \\
                   &=\lambda \lim_{h \to 0^\pm}(\lambda h)^{-1}(\diri(u+h\lambda v)-\diri(u)) = \\
                   &=\lambda \Lambda^\pm(u,v).
\end{align*} 
Let $u,v \in D(\E),$ let $\lambda >0$. Hence
\begin{align*}
\Lambda^\pm(\lambda u,v) &= \lim_{h \to 0^\pm} h^{-1}(\diri(\lambda u+hv)-\diri(\lambda u)) = \\
                   &= \lim_{h \to 0^\pm} h^{-1}(\diri(\lambda(u+h\lambda^{-1}v))-\diri(\lambda u)) = \\
                   &= \lambda \lim_{h \to 0^\pm} \lambda h^{-1}(\diri(u+h\lambda^{-1}v))-\diri(u)) = \\
                   &= \lambda \Lambda^\pm(u,v).
\end{align*}
We prove only the convexity of $\Lambda^+(u,\cdot)$, being the other proof very similar. Fix $u,v_1,v_2 \in D(\E)$.
\begin{align*}
&\Lambda^+(u, \lambda v_1 + (1-\lambda) v_2) = \\
&= \lim_{\sigma \to 0^+} \sigma^{-1}(\diri(u+\sigma \lambda v_1 + \sigma(1-\lambda)v_2)-\diri(u)) =\\
&= \lim_{\sigma \to 0^+} \sigma^{-1}(\diri(\lambda(u+\sigma  v_1) + (1-\lambda)(u+\sigma v_2)-\diri(u)) \leq \\
&\leq \lim_{\sigma \to 0^+} \sigma^{-1}(\lambda \diri(u+\sigma v_1) + (1-\lambda) \diri(u+\sigma v_2) - \diri(u)) = \\
&= \lambda \Lambda^+(u,v_1) + (1-\lambda)\Lambda^+(u,v_2).
\end{align*}
$5.$\\
It is sufficient to switch $h$ with $-h$ and take limits.

\noindent $6.$\\
The proof of the last statement is a direct calculation. We perform it for the two limits at once: 
\begin{align*}
\lim_{h \to 0} \dfrac{\diri(u+hv)-\diri(u)}{h} =  \lim_{h \to 0} \dfrac{\diri(hv)}{h} =  \lim_{h \to 0} \dfrac{h^2\diri(v)}{h}=0.
\end{align*}
    
\end{proof}

Note that $\Lambda^\pm(u,v)$ generalise the integrals $\int_X D^\pm u(\nabla v) \, dm$  introduced by Gigli in  \cite{gigli2015differential} for metric measure spaces. There, even the pointwise objects $D^\pm u(\nabla v)(x)$ make sense, while in our setting the slopes $\Lambda^{\pm}$ are not necessarily represented by a density w.r.t. $dm$. 

In general, $\Lambda^- \neq \Lambda^+$, even in Finsler manifolds. In case $\Lambda^+(u,\cdot) = \Lambda^-(u,\cdot)$, the form $\E$ is said to be \emph{regular} at $u$, and more structure is available. We also say that $\E$ is \emph{regular} if $\E$ is regular at all $u \in D(\E).$ If $\E$ is Fr\'echet-differentiable, we have that $\E$ is regular and $\Lambda(u,v) = (\nabla E(u),v),$ for all $u,v \in D(\E).$

\begin{prop}\label{pro:pro1}
    Let $\E$ be a $2$-homogeneous, regular nonlinear Dirichlet form. 
    Then, $\Lambda := \Lambda^+=\Lambda^-$ is linear in the second argument. Moreover, for all $u \in D(\E)$, we have that $\Lambda(u,\cdot) \in D(\E)^\star$.
    Finally, $\E$ is quadratic if and only if $\E$ is regular and 
    \begin{equation}\label{eq:eq14}
    \forall u,v \in D(\E), \qquad   \Lambda(u,v) = \Lambda(v,u). 
    \end{equation}
\end{prop}

\begin{proof}
    The first assertion is entailed by the fact that $\Lambda(u,\cdot)$ is both concave and convex, and continuous, see Theorem \ref{thm:thm5}. 
    If $\E$ is quadratic, we have 
    $$\Lambda(u,v) = \int_X \sqrt{A}(u) \, \sqrt{A}(v) \, dm,$$
    which is symmetric in $u,v.$
    For the converse, assume $\E$ to be regular. Then, the maps $t \mapsto \E(u+tv)$ are differentiable, for all $u,v \in D(\E).$ Then,
    \begin{align*}
    &\E(u+v) - \E(u) = \int_0^1 \frac{d}{dt} \, \E(u+tv) \, dt = \\
    &= \int_0^1 \Lambda(u+tv,v) \, dt = \int_0^1 \Lambda(v,u) + 2t \, \Lambda(v,v) = \Lambda(v,u) + \E(v). 
    \end{align*}    
    By exchanging $v$ with $-v$, and adding up, we prove that $\E$ satisfies the parallelogram identity, hence $\|\cdot\|_{\E}$ is induced by a scalar product. Equivalently, $\E$ is quadratic.  
\end{proof}

The last point shows that the arguments $u,v$ of $\Lambda$ actually play non-interchangeable roles (unless $\Lambda$ is a bilinear Dirichlet form). This is intuitive in the case of Finsler manifolds, see \eqref{eq:eq12}.

\section{Second-order calculus}\label{sec:sec6}

In metric measure spaces, the role of the Laplacian is played by the minimal section of $\partial \mathrm{Ch}$, being $\mathrm{Ch}$ the Cheeger's energy.
In this section, we investigate the properties of $\partial \E$, for a $2$-homogeneous nonlinear Dirichlet form. Moreover, we introduce an extended subdifferential, in analogy with \cite{gigli2015differential}.

Let $\E$ be a $2-$homogeneous nonlinear Dirichlet form. 
Consider the space $D(\partial \E)$, equipped with the distance 
$$d^2_\partial(u,v) = \|u-v\|^2_{\mathrm{L}^2(X,m)} + \|\partial \E(u) - \partial \E(v)\|^2_{\mathrm{L}^2(X,m)},$$
where the second contribution is a distance between closed subsets of $\mathrm{L}^2(X,m).$

\begin{prop}\label{prop:prop2}
    Let $\E$ be a $2-$homogeneous nonlinear Dirichlet form. Then, we have that $D(\partial \E) \subset D(\E)$ is dense and continuous at $0$. Moreover, $D(\partial \E) \subset \mathrm{L}^2(X,m)$ is dense and continuous. Finally, 
    $$\forall u \in D(\partial \E), \qquad \E(u) \leq \frac{1}{2} \, \|\partial^0 \E( u)\|_{\mathrm{L}^2(X,v)}  \, \|u\|_{\mathrm{L}^2(X,m)}.$$
\end{prop}

\begin{proof}
    The density of the two inclusions is proved as follows. Let $u \in \mathrm{L}^2(X,m)$. Then, $T_t u \in D(\partial \E)$ for all $t>0$ \cite{evans1998partial}, where $(T_t)_t$ is the gradient flow generated by $\E.$ The convergence properties of $T_t u \to u,$ as $t\to 0$ are sufficient to conclude. 
    In addition, $d_\partial(u,v) \geq \|u-v\|_{\mathrm{L}^2(X,m)}$ shows the continuity of the inclusion $D(\partial \E) \subset \mathrm{L}^2(X,m)$. 
    The inequality $\E(u) \leq \frac{1}{2} \, \|\partial^0 \E( u)\|_{\mathrm{L}^2(X,v)}  \, \|u\|_{\mathrm{L}^2(X,m)},$ is a combination of \eqref{cdc2} and Cauchy-Schwarz's inequality. This implies also the continuity at $0$ of $D(\partial \E) \subset D(\E).$
\end{proof}

Some integration by parts rules, reminiscent of those given in \cite{ambrosio2014calculus,gigli2015differential}, link the subdifferential with the slopes of Section \ref{sec:sec5}, as the next result shows.   

\begin{thm}
    \label{thm:thm6}
    Let $\E$ be a $2$-homogenous nonlinear Dirichlet form. Then, the following rules hold. 
    \begin{enumerate}
        \item For all $u \in D(\partial E),$ and all $v \in D(\E),$ we have that
        $$ \Lambda^-(u,v) \leq (\partial \E(u),v) \leq \Lambda^+(u,v).$$
        In particular, if $\E$ is such that $\Lambda(u,\cdot)^+ = \Lambda^-(u,\cdot),$ and $\E$ is subdifferentiable at $u$, we have 
        $$\forall v \in D(\E), \qquad \Lambda^\pm(u,v) = (\partial \E(u),v),$$
        and $\partial \E(u)$ contains only one element.
        \item For $\lambda > 0$, let $A_\lambda$ be the Yosida regularisation of $\partial \E.$ Then, 
        $\forall u,v \in D(\E)$, 
        $$\Lambda^-(u,v) \leq \liminf_{\lambda \to 0} (A_\lambda (u),v) \leq \limsup_{\lambda \to 0 } (A_\lambda (u),v) \leq \Lambda^+(u,v).$$
        In particular, if $\E$ is regular, we have $ \lim_{\lambda \to 0} (A_\lambda(u),v) = \Lambda(u,v),$ for all $u,v \in D(\E).$
    \end{enumerate}
\end{thm}

\begin{proof}
    $1.$\\
We perform the calculation for one side of the inequality, being the other one analogous.
Fix $\xi \in \partial \diri(u)$ and compute
$$\diri(u+hv)-\diri(u) \geq (\xi,hv),$$
which reads, if $h<0$, as
$$h^{-1}(\diri(u+hv)-\diri(u)) \leq (\xi,v).$$
The inequality follows by taking the supremum for $h<0$.
If $\E$ is regular at $u,$ we have that $(\partial \E(u),v)$ is prescribed for all $v \in D(\E)$ by the values of $\Lambda(u,v)$, hence, $\partial \E(u)$ contains only one element.

\noindent $2.$\\

We start with the $\liminf$ inequality. 
For all $\lambda > 0$, $u,v$ as in the hypothesis and $h<0$, consider
$$\diri_\lambda (u+hv) - \diri_\lambda(u) \geq (A_\lambda(u),hv),$$
which leads to
$$\liminf_{\lambda \to 0} h^{-1} (\diri_\lambda (u+hv) - \diri_\lambda(u)) \leq \liminf_{\lambda \to 0} (A_\lambda(u),v).$$
The l.h.s. admits a limit, hence,
$$h^{-1} (\diri(u+hv)-\diri(u)) \leq \liminf_{\lambda \to 0} (A_\lambda(u),v).$$
Taking the supremum over $h<0$ yields the sought inequality.
For the $\limsup$ inequality, still consider $\lambda>0$, $u,v$ as in the hypotheses, and $h>0$. 
Write 
$$\diri_\lambda(u+hv)-\diri_\lambda(u) \geq (A_\lambda u,hv),$$
and take the $\limsup$ in both sides to get
$$\limsup_{\lambda \to 0} \diri_\lambda(u+hv)-\diri_\lambda(u) \geq \limsup_{\lambda \to 0} (A_\lambda u,hv),$$
which reads 
$$\diri(u+hv)-\diri(u) \geq  \limsup_{\lambda \to 0} (A_\lambda u,hv),$$
hence, 
$$h^{-1} \diri(u+hv)-\diri(u) \geq  \limsup_{\lambda \to 0} (A_\lambda u,v).$$
The result follows by taking the infimum in $h$.
\end{proof}

The domain of $\partial \E$ generalises the space $\mathrm{W}^{2,2}(\mathbb R^d)$, with $$\partial \E : D(\partial \E) \to \mathrm{L}^2(X,m).$$
Hereby, we give an extended definition, in the spirit of \cite{gigli2015differential}, mimicking the distributional Laplacian $$\Delta: \mathrm{W}^{1,2}(\mathbb R^d) \to \mathrm{H}^{-1}(\mathbb R^d).$$

Let $\E$ be a $2$-homogeneous nonlinear Dirichlet form. 
Then, we say that a function $u \in D(\E)$ is a point of extended subdifferentiability if there exists a measure $\mu$ on $(X,\mathcal{F})$ such that all $v \in D(\E)$ are $\mu$-measurable and the following holds
\begin{equation}
    \label{eq:eq15}
    \Lambda^{-}(u,v) \leq \int_X v \, d\mu \leq \Lambda^+(u,v).
\end{equation}
In this case, we write $u \in D(\bar{\partial} \E)$ and  $\mu \in \bar{\partial} \E(u)$. Notice that Theorem \ref{thm:thm5} implies $\mu \in D(\E)^\star.$

We collect some properties of $\bar{\partial} \E$ in the next result, which concludes our analysis. 

\begin{prop}\label{thm:thm7}
    Let $\E$ be a $2$-homogeneous nonlinear Dirichlet form. Then, the extended subdifferential $\bar{\partial}$ satisfies the following:
    \begin{enumerate}
        \item $\bar{\partial} \E(u)$ is convex and $1$-homogeneous;
        \item if $u \in D(\bar{\partial}\E)$ and $\E$ is regular at $u$, then $\bar{\partial}\E(u)$ contains only one measure.
        \item if $u\in D(\partial \E),$ and $\xi \in \partial \E(u)$, then $\xi \, dm \in \bar{\partial}\E(u).$ 
    \end{enumerate}
\end{prop}

\begin{proof}
    The first assertion is a consequence of the $1$-homogeneity of $\Lambda^\pm$ and of the linearity of the integral with respect to the measure. 
    The second assertion follows from the fact that the values of $\int_X v \, d\mu$ are prescribed by $\Lambda(u,v).$ Finally, the third statement holds by definition, after Theorem \ref{thm:thm6}. 
\end{proof}

The validity of a converse of the third statement in the last theorem has been discussed for metric spaces  in \cite{gigli2015differential}, and looks unclear. We are not able to give an analogous of \cite[Proposition 4.11]{gigli2015differential} as well. Such formula is one of the hypotheses of \cite{ambrosio2015bakry} for the quadratic case.  

\section{Perspectives and \emph{desiderata}}\label{sec:sec7}

A missing point in the theory is a pointwise object representing $\Lambda^{\pm}(u,v),$ which should play the same role as $D^{\pm} v (\nabla u)(x)$ in \cite{gigli2015differential}. It is unclear whether the existence of a density for $\Lambda^{\pm}$ should be imposed, or if it is a consequence of some locality hypotheses on $\E$. In any case, heuristically, it corresponds to finer integration by parts formulae as those of the current paper. Also a Leibniz formula and some chain rules on $\Lambda^\pm(u,\cdot)$ would be advances in the theory. Still, we do not know if one should expect such properties or impose them. 

Another desirable statement is some integral representation of $\E$ under locality assumptions. However, this is not available even if $X = \mathbb R^2.$ Finally, a missing notion is that of \emph{linearised} flow associated to the gradient flow of a nonlinear Dirichlet form. Equivalently, one would need a \emph{linearisation of $\partial \E$ in the direction of the gradient}, as in \cite{ohta2014bochner}.

\let\d\oldd
\let\k\oldk
\let\l\oldl
\let\L\oldL
\let\o\oldo
\let\O\oldO
\let\r\oldr
\let\S\oldS
\let\t\oldt
\let\u\oldu

\section*{Declarations}
The author has been funded by the European Union’s Horizon 2020 research and innovation program under the Marie Skłodowska-Curie grant agreement No 754362. Partial support has been obtained from the EFI ANR-17-CE40-0030 Project of the French National Research Agency. 
\\
Data sharing not applicable to this article as no datasets were generated or analysed during the current study.
\\
The author has no competing interests to declare that are relevant to the content of this article.
\\
The author thanks the organisers of the $4^{th}$ BYMAT Conference for the opportunity.

\bibliography{ref.bib}
\bibliographystyle{siam}
\end{document}